# Solving equality-constrained optimization problems without Lagrange multipliers

Cyril Cayron[1]

[1] Ecole Polytechnique Fédérale de Lausanne (EPFL), Laboratoire de métallurgie thermomécanique (LMTM), PX-Group chair

Email: cyril.cayron@epfl.ch    ORCID id: 0000-0002-3190-5158

Abstract: Constrained optimization problems are classically solved with the help of the Lagrange multipliers and the Lagrangian function. However, the disadvantage of this approach is that it artificially increases the dimensionality of the problem. Here, we show that the determinant of the Jacobian of the problem (function to optimize and constraints) is null. This extra equation directly transforms any equality-constrained optimization problem into a solving problem of same dimension. We also introduced the constraint matrices as the largest square submatrices of the Jacobian of the constraints. The boundaries of the constraint domain are given by the nullity of their determinants. The constraint matrices also permit to write the function to be optimized as a Taylor series of any of its variable, with its coefficients algebraically determined by an iterative process of partial derivation.

## 1    Introduction

Constrained optimization problems exist nearly everywhere in science. In mechanics, the trajectory of a particle is that one that optimizes the action (the difference between the kinetics energy and the potential energy) along the path with the constraint that the total energy (the sum of the kinetics energy and the potential energy) is constant [1,2]. In thermodynamics, the engineer wants to extract the maximum external work of a system to run a mechanical machine. The optimum is obtained by reducing at its minimum the irreversible entropy created during the heat exchanges between the system and the machine with the constraint that the sum of system's and machine's energies is constant [3,4]. Inversely, if the system is free to spontaneously change, no work can be extracted and the system's entropy increases up to its maximum. In economics, one would like a system in which the exchanges of goods, information etc. are such that they maximize the effective work useful for the society while minimizing the financial transactions with their inevitable losses in tax havens, with the constraints that the total amount of goods is limited and the amount of money is controlled by the central banks. In general these problems are very complex because of the high number of parameters [5,6]. During a pandemic crisis, such as that of Coronavirus we are facing at the moment, trying to maintain at the minimum the mortality while keeping the economy running is also an optimization problem with a high number of intricate parameters.

The classical method to solve constrained optimization problems consists in using the Lagrange multipliers [7-9]. Lagrange multiplier theorem states that for a real function $f$ in a $N$-dimension space of parameters $x_1, \ldots, x_N$ that follows $M$ constrains $g_k(x_1, \ldots, x_N) = C_k$ for $k \in [1, M]$, the stationary points (maximum, minima or inflexion) are such that in these points, $\nabla f$ the gradient of $f$ is a linear combination of $\nabla g_k$ the gradients of the constrained function $g_k$. The Lagrange multipliers $\lambda_k$ are the coefficients of this linear relation:

$$\nabla f(x_1, \ldots, x_N) = \sum_{k=1}^{M} \lambda_k \nabla g_k(x_1, \ldots, x_N) \tag{1}$$

where the gradients $\nabla f$ and $\nabla g_k$ are $N$-dimension vectors.

Let us illustrate the Lagrange theorem in 2D with a function $f(x, y)$. The function $f$ can represented by contour curves $f(x, y) = c$. At the optimum points obtained for some specific values of $c$, the curve becomes tangent to that of the constraint $g(x, y) = C$, which implies than the gradient vectors $\nabla f$ and $\nabla g$ are parallel. Two graphical examples are shown in Fig. 1 with

(1a) $f(x, y) = x^2 + 2y^2$ with the constraint $g(x, y) = x - y^2 = 1$

(1b) $f(x, y) = x^2 + 2y^2 - 2xy^2$ with the constraint $g(x, y) = x + y^2 = 1$

The contour curves of $f(x, y)$ are represented in blue-to-red thermal colors, and the contour curve of the constraint $g(x, y)$ in red. In the example (1a), there is one solution $(x_1, y_1) = (1,0)$ for which $f(x_1, y_1) = 1$ enlightened by the blue curve. In the example (1b), there are three solutions $(x_1, y_1) = (1,0)$, $(x_2, y_2) = \left(\frac{2}{3}, \frac{1}{\sqrt{3}}\right)$ and $(x_3, y_3) = \left(\frac{2}{3}, \frac{-1}{\sqrt{3}}\right)$ for which the values of $f(x, y)$ are respectively 1 (curve in cyan), twice 2/3 (curves in blue). The details of the calculation will be given in section 4.

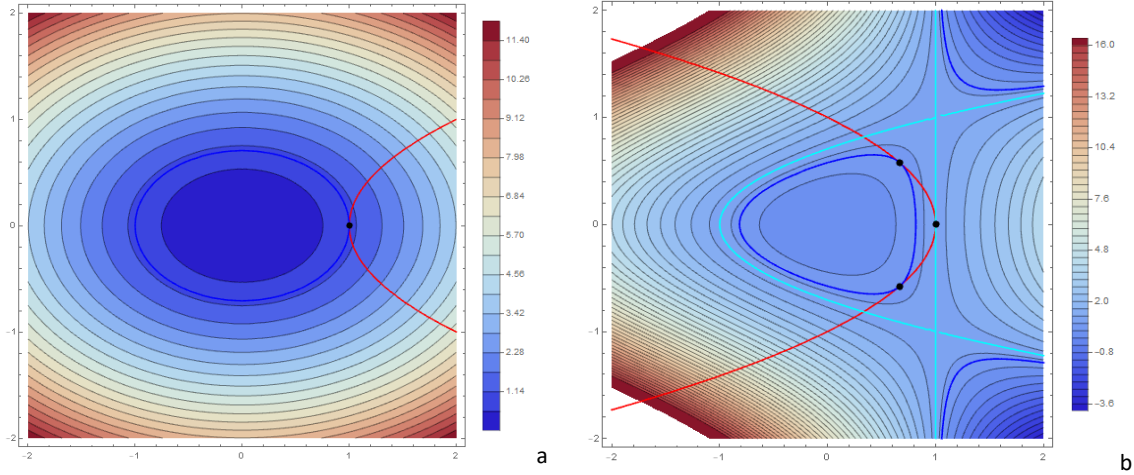

Fig. 1. Two 2D examples representing the function $f(x, y)$ by its contour curves $f(x, y) = c$. The optimum points are obtained for specifc values of $c$ when the contour curv of $f$ becomes tangent to the contour curve of $g(x, y) = C$. At these points, marked by the black spots, the gradients $\nabla f$ and $\nabla g$ are parallel and linearly linked by the Lagrange multipliers. (a) $f(x, y) = x^2 + 2y^2$ with the constraint $g(x, y) = x - y^2 = 1$, (b) $f(x, y) = x^2 + 2y^2 - 2xy^2$ with the constraint $g(x, y) = x + y^2 = 1$.

The usual optimization method consists in introducing the Lagrangian function $\mathcal{L}$

$$\mathcal{L}(x_1, \ldots, x_N, \lambda_1, \ldots, \lambda_M) = f(x_1, \ldots, x_N) - \sum_{k=1}^{M} \lambda_k g_k(x_1, \ldots, x_N) \tag{2}$$

The optimum points of $f$ are the solutions of a system of $N + M$ equations with $N + M$ unknown $\nabla \mathcal{L}(x_1, \ldots, x_N, \lambda_1, \ldots, \lambda_M) = 0$ and $\forall k \in [1, M], g_k(x_1, \ldots, x_N) = C_k$, or more explicitly

$$\begin{cases} \dfrac{\partial f}{\partial x_1}(x_1,\dots,x_N) - \sum_{k=1}^{M} \lambda_k \dfrac{\partial g_k}{\partial x_1}(x_1,\dots,x_N) = 0 \\ \qquad\vdots \\ \dfrac{\partial f}{\partial x_N}(x_1,\dots,x_N) - \sum_{k=1}^{M} \lambda_k \dfrac{\partial g_k}{\partial x_N}(x_1,\dots,x_N) = 0 \\ g_1(x_1,\dots,x_N) = C_1 \\ \qquad\vdots \\ g_M(x_1,\dots,x_N) = C_M \end{cases} \qquad (3)$$

Consequently, the use of the Lagrange multipliers in optimization problems "artificially" increases the space dimension from $N$ to $N + M$. In addition, with this method, the number of constrains $M$ seems to be decorrelated from the dimensionality of the problem $N$. This is unfortunate because we know that we need $N - 1$ constraint equations that are non-linearly linked to find the points of optimization. The method with Lagrange multipliers leads to solve 2$N$-1 equations. We will show that the problem can be solved in the initial $N$-dimension space, without the Lagrange multipliers and without increasing the dimensionality of the problem. The method we propose results from an intuitive, simple and effective equation. All the graphics that will be shown were plotted with Mathematica 12.

## 2  Making the optimization problem a solving problem of same dimension

### 2.1  Intuitive approach

Let us consider a 2D problem. We want to optimize $f(x, y)$ with the constraint $g(x, y) = C$. Instead of representing the problem in 2D we temporarily use a third dimension $\varphi$ to plot $\varphi = f(x, y)$. In this space the constraint $g(x, y) = C$ appears as a vertical cylinder along $\varphi$. The points with $\varphi = f(x, y)$ and with $g(x, y) = C$ form two surfaces that intersect into a curve in the $(x, y, \varphi)$ space, as illustrated in the examples (1a) and (1b) in Fig. 2.

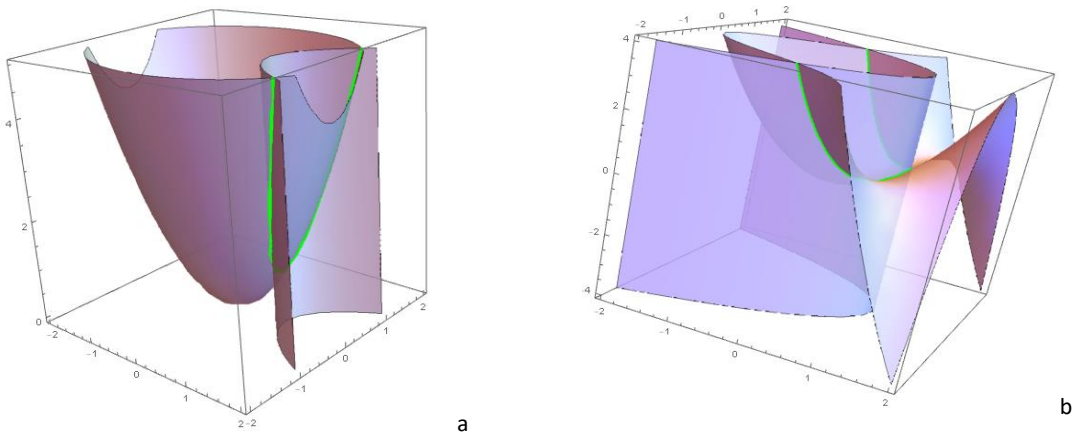

*Fig. 2.  3D surface representation the 2D problems previously shown by iso-contours in Fig. 1. In the $(x, y, \varphi)$ space (here $\varphi$ is te vertical direction), the intersection of the surfaces represeting the function $\varphi = f(x, y)$ with the cylinder representing constraint $g(x, y) = C$ is a curve enlightned in green. Its extrema are the solution of the constrained optimum problem.*

The optimum points of this 3D curve are the optimum points of $f(x, y)$. They are given by $d\varphi = \dfrac{\partial f(x,y)}{\partial x} dx + \dfrac{\partial f(x,y)}{\partial y} dy = 0$. As we are working now in the derivative space, let us express the constraint condition by its derivative $\frac{\partial g(x,y)}{\partial x}dx + \frac{\partial g(x,y)}{\partial y}dy = 0$. The optimum and constraint conditions can be written in a matrix form as

$$\begin{bmatrix} \frac{\partial f(x,y)}{\partial x} & \frac{\partial f(x,y)}{\partial y} \\ \frac{\partial g(x,y)}{\partial y} & \frac{\partial g(x,y)}{\partial y} \end{bmatrix} \begin{bmatrix} dx \\ dy \end{bmatrix} = \begin{bmatrix} 0 \\ 0 \end{bmatrix}$$

with a non-null infinitesimal vector $(dx, dy)$. Consequently

$$Det \begin{bmatrix} \frac{\partial f(x,y)}{\partial x} & \frac{\partial f(x,y)}{\partial y} \\ \frac{\partial g(x,y)}{\partial y} & \frac{\partial g(x,y)}{\partial y} \end{bmatrix} = 0 \qquad (4)$$

In other words, the optimum condition is replaced by a "constraint-like" equation, without Lagrange multiplier, without increasing the dimensionality of the initial problem.

## 2.2 Generalization

The generalization is direct. Let us consider a space of dimension $N$ with a real function $f(x_1, ..., x_N)$ that we would like to optimize while respecting the $N-1$ constraints $g_k(x_1, ..., x_N) = C_k$ for $k \in [1, N-1]$. We use $(f, g_1, ... g_{N-1})$ the $N$-vector constituted by the functions $f$ and all the $g_k$ to define the Jacobian of the problem $J_{(f,g)}$ by

$$J_{(f,g)} = \begin{bmatrix} \frac{\partial f(x_1, ..., x_N)}{\partial x_1} & \cdots & \frac{\partial f(x_1, ..., x_N)}{\partial x_N} \\ \frac{\partial g_1(x_1, ..., x_N)}{\partial x_1} & \cdots & \frac{\partial g_1(x_1, ..., x_N)}{\partial x_N} \\ \vdots & \ddots & \vdots \\ \frac{\partial g_{N-1}(x_1, ..., x_N)}{\partial x_1} & \cdots & \frac{\partial g_{N-1}(x_1, ..., x_N)}{\partial x_N} \end{bmatrix} \qquad (5)$$

The points $(x_1, ..., x_N)$ on which $f$ is optimum with the constraints $g_k$ are those that check

$$J_{(f,g)} \begin{bmatrix} dx_1 \\ \vdots \\ dx_N \end{bmatrix} = \begin{bmatrix} 0 \\ \vdots \\ 0 \end{bmatrix}$$

with a non-null infinitesimal vector $(dx_1, ..., dx_N)$. Consequently, the "constraint-like" equation that substitutes the optimum condition on $f$ to complete the $N-1$ initial constraints $g_k$ is simply

$$Det\, J_{(f,g)} = 0 \qquad (6)$$

One can note that equation $Det\, J_{(f,g)} = 0$ leads to an additional equation if and only if the $N-1$ initial constraint equations are not linearly linked. It is also clear that if the number of constraint equations were strictly lower than $N-1$, the problem could not be solved because of its under-determination; the solutions could not be expressed as a finite set of points.

Now, let us consider again the Lagrange theorem that tells that the optima are obtained when the gradient of $f$ is a linear combination of the gradients of the constraints $g_k$, i.e. $\nabla f(x_1, ..., x_N) = \sum_{k=1}^{M} \lambda_k \nabla g_k(x_1, ..., x_N)$. Writing the gradients as line vectors, it is clear that such linearity condition is equivalent to the condition $Det(J_F) = 0$. This condition makes any optimization problem of dimension

$N$ equivalent to a solving a problem of $N$ equations. With the classical method based on the Lagrange multipliers and its associated Lagrangian function, one would have to solve $2N-1$ equations. Equation (6) is directly the result of the additional N-1 equations introduced with the Lagrange multipliers. The fact that Lagrange multipliers are not required for solving optimization problems was already proposed by Gigena by using the implicit function theorem [10].

The system of $N$ equations (5) now allows us to determine the optima of the function $f$. If the expression of $f$ and $g$ are complex, the solutions can be obtained only numerically. However, we will show that the function $f$ can be developed in Taylor series of any of the variable $x_k$ with coefficients algebraically determined thanks the constraint matrices. These matrices are also useful to determine the boundaries of the variables $(x_1, ..., x_N)$ imposed by the constraints $g$.

## 3  The constraint matrices, the *g*-constraint domain and the Taylor series of *f*

### 3.1  The constraint matrices and the infinitesimal variables

Let us consider here only the constraints. The infinitesimal entities $dx_k$ are linearly linked each other by the derivative of the constraint equations

$$\begin{bmatrix} \frac{\partial g_1(x_1,...,x_N)}{\partial x_1} & \cdots & \frac{\partial g_1(x_1,...,x_N)}{\partial x_N} \\ \vdots & \ddots & \vdots \\ \frac{\partial g_{N-1}(x_1,...,x_N)}{\partial x_1} & \cdots & \frac{\partial g_{N-1}(x_1,...,x_N)}{\partial x_N} \end{bmatrix} \begin{bmatrix} dx_1 \\ \vdots \\ dx_N \end{bmatrix} = \begin{bmatrix} 0 \\ \vdots \\ 0 \end{bmatrix} \quad (7)$$

The $(N-1)N$ matrix in this equation is the Jacobian $J_{(g)}$ of the vectorial function $(g_1, ... g_{N-1})$, the system is 1-underdetermined, which means that any $dx_k$ is a linear combination of the other $dx_i$ with $i \neq k$. The practical way to establish the linear relation is simple. Let us explain it with $dx_1$. Equation (7) can be written

$$\begin{bmatrix} \frac{\partial g_1(x_1,...,x_N)}{\partial x_1} \\ \vdots \\ \frac{\partial g_{N-1}(x_1,...,x_N)}{\partial x_1} \end{bmatrix} dx_1 + \begin{bmatrix} \frac{\partial g_1(x_1,...,x_N)}{\partial x_2} & \cdots & \frac{\partial g_1(x_1,...,x_N)}{\partial x_N} \\ \vdots & \ddots & \vdots \\ \frac{\partial g_{N-1}(x_1,...,x_N)}{\partial x_2} & \cdots & \frac{\partial g_{N-1}(x_1,...,x_N)}{\partial x_N} \end{bmatrix} \begin{bmatrix} dx_2 \\ \vdots \\ dx_N \end{bmatrix} = \begin{bmatrix} 0 \\ \vdots \\ 0 \end{bmatrix}$$

Let us call $S_1(x_1, ..., x_N) = \begin{bmatrix} \frac{\partial g_1(x_1,...,x_N)}{\partial x_2} & \cdots & \frac{\partial g_1(x_1,...,x_N)}{\partial x_N} \\ \vdots & \ddots & \vdots \\ \frac{\partial g_{N-1}(x_1,...,x_N)}{\partial x_2} & \cdots & \frac{\partial g_{N-1}(x_1,...,x_N)}{\partial x_N} \end{bmatrix}$.

It is a $(N-1)^2$ square submatrix of the Jacobian $(N-1)N$ matrix $J_{(g)}$. If at the optimum point $(x_1, ..., x_N)$, the matrix $S_1$ is invertible, then all the other infinitesimals $dx_i$ with $i \neq 1$ are proportional to $dx_1$ by the set of equalities written in the vectorial form

$$\begin{bmatrix} dx_2 \\ \vdots \\ dx_N \end{bmatrix} = -S_1^{-1}(x_1, ..., x_N) \begin{bmatrix} \frac{\partial g_1(x_1,...,x_N)}{\partial x_1} \\ \vdots \\ \frac{\partial g_{N-1}(x_1,...,x_N)}{\partial x_1} \end{bmatrix} dx_1 = \begin{bmatrix} s_{2,1}(x_1,...,x_N) \, dx_1 \\ \vdots \\ s_{N,1}(x_1,...,x_N) \, dx_1 \end{bmatrix}$$

The same reasoning can be applied to any infinitesimal $dx_k$. We call $S_k$ the submatrix of $J_{(g)}$ obtained by removing the *k*-line and the *k*-column:

$$S_k(x_1,\ldots,x_N) = \begin{bmatrix} \dfrac{\partial g_1}{\partial x_1} & \cdots & \dfrac{\partial g_1}{\partial x_{k-1}} & \dfrac{\partial g_1}{\partial x_k} & \dfrac{\partial g_1}{\partial x_{k+1}} & \cdots & \dfrac{\partial g_1}{\partial x_N} \\ \vdots & \vdots & \vdots & \vdots & \vdots & \vdots & \vdots \\ \dfrac{\partial g_{k-1}}{\partial x_1} & \cdots & \dfrac{\partial g_{k-1}}{\partial x_{k-1}} & \dfrac{\partial g_{k-1}}{\partial x_k} & \dfrac{\partial g_{k-1}}{\partial x_{k+1}} & \cdots & \dfrac{\partial g_{k-1}}{\partial x_N} \\ \dfrac{\partial g_k}{\partial x_1} & \cdots & \dfrac{\partial g_k}{\partial x_{k-1}} & \dfrac{\partial g_k}{\partial x_k} & \dfrac{\partial g_k}{\partial x_{k+1}} & \cdots & \dfrac{\partial g_k}{\partial x_N} \\ \dfrac{\partial g_{k+1}}{\partial x_1} & \cdots & \dfrac{\partial g_{k+1}}{\partial x_{k-1}} & \dfrac{\partial g_{k+1}}{\partial x_k} & \dfrac{\partial g_{k+1}}{\partial x_{k+1}} & \cdots & \dfrac{\partial g_{k+1}}{\partial x_N} \\ \vdots & \vdots & \vdots & \vdots & \vdots & \vdots & \vdots \\ \dfrac{\partial g_{N-1}}{\partial x_1} & \cdots & \dfrac{\partial g_{N-1}}{\partial x_{k-1}} & \dfrac{\partial g_{N-1}}{\partial x_k} & \dfrac{\partial g_{N-1}}{\partial x_{k+1}} & \cdots & \dfrac{\partial g_{N-1}}{\partial x_N} \end{bmatrix}$$

Equation (7) becomes

$$\begin{bmatrix} \dfrac{\partial g_1(x_1,\ldots,x_N)}{\partial x_k} \\ \vdots \\ \dfrac{\partial g_{N-1}(x_1,\ldots,x_N)}{\partial x_k} \end{bmatrix} dx_k + S_k \begin{bmatrix} dx_1 \\ \vdots \\ dx_{k-1} \\ dx_k \\ dx_{k+1} \\ \vdots \\ dx_N \end{bmatrix} = \begin{bmatrix} 0 \\ \vdots \\ 0 \end{bmatrix} \quad (8)$$

If $S_k$ is invertible, the infinitesimals $dx_i\ \forall i \neq k$ are proportional to $dx_k$ by

$$\begin{bmatrix} dx_1 \\ \vdots \\ dx_{k-1} \\ dx_k \\ dx_{k+1} \\ \vdots \\ dx_N \end{bmatrix} = -S_k^{-1} \begin{bmatrix} \dfrac{\partial g_1}{\partial x_k} \\ \vdots \\ \dfrac{\partial g_{N-1}}{\partial x_k} \end{bmatrix} dx_k = \begin{bmatrix} s_{1,k}\, dx_k \\ \vdots \\ s_{k-1,k}\, dx_k \\ s_{k,k}(x_1,\ldots,x_N)\, dx_k \\ s_{k+1,k}\, dx_k \\ \vdots \\ s_{N,k}\, dx_k \end{bmatrix} \quad (9)$$

Consequently, for any fixed $k$, $\forall i \neq k$,

$$dx_i = s_{i,k} dx_k \quad \text{with } s_{i,k} = -(S_k^{-1})_{i,l} \dfrac{\partial g_l(x_1,\ldots,x_N)}{\partial x_k} \quad (10)$$

The equation was written here in a compact form, but $s_{i,k}$ are actually functions of the $N$ variables, and one should read for them $s_{i,k}(x_1,\ldots,x_N)$.

### 3.2   The constraint matrices and the extrema of the g-constraint curve

The set of the $N-1$ constraint equations $g_k(x_1,\ldots,x_N) = C_k$ define the constraint curve $\mathcal{C}_{(g)}$. This curve is bounded in the $(x_1,\ldots,x_N)$ space. For any fixed $k$, the extrema values on the $x_k$-axis verify the condition $dx_k = 0$, with $dx_i \neq 0$ for at least one $i \neq k$. If the vector $\begin{bmatrix} \dfrac{\partial g_1(x_1,\ldots,x_N)}{\partial x_k} \\ \vdots \\ \dfrac{\partial g_{N-1}(x_1,\ldots,x_N)}{\partial x_k} \end{bmatrix}$ is also not null, the condition (8) becomes $Det\ S_k = 0$. Consequently, the extrema of the constraint curve $\mathcal{C}_{(g)}$ along the $x_k$-axis are the points $(x_1,\ldots,x_N)$ that are the solutions of the system of $N-1$ equations

$$\begin{cases} g_1(x_1, \ldots, x_N) = C_1 \\ \quad \vdots \\ g_{k-1}(x_1, \ldots, x_N) = C_{k-1} \\ \text{Det } S_k(x_1, \ldots, x_N) = 0 \\ g_{k+1}(x_1, \ldots, x_N) = C_{k+1} \\ \quad \vdots \\ g_{N-1}(x_1, \ldots, x_N) = C_{N-1} \end{cases} \tag{11}$$

## 3.3 The constraint matrices and the Taylor's series of $f$

The Lagrange multipliers are the partial derivatives of the Lagrangian on the constraint value of $g_k(x_1, \ldots, x_N) = c_k$, i.e. $\lambda_k = \frac{\partial \mathcal{L}(x_1, \ldots, x_N, \lambda_1, \ldots, \lambda_M)}{\partial c_k}$. They are thus often used to evaluate the rate of change of the optima due a relaxation of a given constraint $c_k$. However, it is also possible to describe the "nature" of optimum due a relaxation of the initial parameters $x_1, \ldots, x_N$; and this can be done without the Lagrange multipliers, by using the constraint matrices that links the infinitesimals. As the constraint equations are not linearly linked, there is at least one matrix $S_k$ that is invertible, which means that the nature of the optimality of $f$ can be unambiguously defined at least along one direction $x_k$. We will show that the function $f$ can be developed in a Taylor series of $x_k$ at any optimum point, and more generally at any point of the constraint curve $\mathcal{C}_{(g)}$. To simplify the rest of the section, let us consider the case where $\det(S_k) \neq 0, \forall k \in [1, N]$.

Using the infinitesimals in mathematics is however always "risky". For example, we tried to use the linear relations of the infinitesimals to calculate the second derivative $d^2 f(x_1, \ldots, x_N) = \sum_{i=1}^{N} \sum_{j=1}^{N} \frac{\partial^2 f(x_1, \ldots, x_N)}{\partial x_i \partial x_j} dx_i dx_j$, by changing $dx_i dx_j$ by the product of $dx_i$ by $dx_j$. We realized on the simple examples (1a) and (1b) that such a naïve approach leads to incorrect results. We interpret this failure by the fact that the notations $dx_i$ in the first derivative and $dx_i dx_j$ in the second derivative are misleading. One should actually understand them as new variables, and one should write $df(x_1, \ldots, x_N, \ldots dx_i \ldots)$, i.e. $df$ as a function a $2N$ variables, and $d^2 f(x_1, \ldots, x_N, \ldots dx_i dx_j \ldots)$, i.e. $d^2 f$ as a function of function of $N + N^2$ variables. Tensor notations $dx_i \otimes dx_j$ would be more appropriate to describe the variables of multi-derivative, and more generally to describe the differentials, as clearly pointed in Ref. [11].

It is however possible to use the linear relations between the infinitesimals to determine the value of $df$, $d^2 f, \ldots d^n f, \ldots$ and consequently establish the Taylor's series of $f$ along any variable $x_k$ at the $f$-optimum points determined in section 2, and more generally at any point of the constraint curve $\mathcal{C}_{(g)}$. Let us come back to the usual partial derivative formula

$$df(x_1, \ldots, x_N, dx_1 \ldots, dx_N) = \sum_{j=1}^{N} \frac{\partial f(x_1, \ldots, x_N)}{\partial x_i} dx_i \tag{12}$$

For any fixed $k \in [1, N]$, we use the relation (10) to write $df$ as a function of $x_1, \ldots, x_N$ and $dx_k$

$$df(x_1, \ldots, x_N, dx_k) = \left( \sum_{i=1}^{N} \frac{\partial f(x_1, \ldots, x_N)}{\partial x_i} s_{i,k}(x_1, \ldots, x_N) \right) dx_k$$

which means that now we can get the derivative of $f$ along $x_k$ as a function of $x_1, \ldots, x_N$:

$$f'_k(x_1, \ldots, x_N) = \frac{df}{dx_k}(x_1, \ldots, x_N) = \sum_{\substack{i=1 \\ i \neq k}}^{N} \frac{\partial f(x_1, \ldots, x_N)}{\partial x_i} s_{i,k}(x_1, \ldots, x_N) \tag{13}$$

We can repeat this derivative process by applying the partial derivative formula to $f'_k$

$$df'_k(x_1, \ldots, x_N, dx_1 \ldots, dx_N) = \sum_{\substack{i=1 \\ i \neq k}}^{N} \frac{\partial f'_k(x_1, \ldots, x_N)}{\partial x_i} dx_i$$

With relation (10), we get $df'_k(x_1, \ldots, x_N, dx_1 \ldots, dx_N) = \left( \sum_{i=1}^{N} \frac{\partial f'_k(x_1,\ldots,x_N)}{\partial x_i} s_{i,k}(x_1, \ldots, x_N) \right) dx_k$. Thus

$$f''_k(x_1, \ldots, x_N) = \frac{df'_k}{dx_k}(x_1, \ldots, x_N) = \sum_{\substack{i=1 \\ i \neq k}}^{N} \frac{\partial f'_k(x_1, \ldots, x_N)}{\partial x_i} s_{i,k}(x_1, \ldots, x_N)$$

With equation (13), it comes

$$f''_k(x_1, \ldots, x_N) = \sum_{\substack{i=1 \\ i \neq k}}^{N} \sum_{\substack{j=1 \\ j \neq k}}^{N} \frac{\partial f(x_1, \ldots, x_N)}{\partial x_i \partial x_j} s_{j,k}(x_1, \ldots, x_N) \, s_{i,k}(x_1, \ldots, x_N) \tag{14}$$

The same process could be repeated $n$-times to calculate the derivative of $f$ at the order $d^{(n)} f_k(x_1, \ldots, x_N)$.

Here, we limit the development to the second order. At any point $P$ $(x_1^P, \ldots, x_N^P)$, the behavior of $f$ along the $x_k$–axis with the other variables $x_i$ obeying the constraints $g_i(x_1, \ldots, x_N) = C_i \, \forall i \neq k$, i.e. $\tilde{f}(x_k) = f(x_k)_{/\, g_i(x_1,\ldots,x_N)=C_i}$ can be approximate by the parabolic form

$$\tilde{f}(x_k) \approx f(x_1^P, \ldots, x_N^P) + f'_k(x_1^P, \ldots, x_N^P)(x_k - x_k^P) + \frac{1}{2} f''_k(x_1^P, \ldots, x_N^P)(x_k - x_k^P)^2 \tag{15}$$

The calculation of the other coefficients of the Taylor series $d^{(n)} f_k(x_1^P, \ldots, x_N^P)$ could possible by repeating the operation. Equation (15) can be useful to evaluate the behavior of the function $f$ around its optimum points. It is important to note that the algebraic expression $d^{(n)} f_k(x_1^P, \ldots, x_N^P)$ can be determined by repeating the process of derivation and variable change described above, which means that the Taylor series of $f$ as a function of $x_k$ can be determined algebraically to any order $n$, even if the forms of $f(x_1, \ldots, x_N)$ and $g_i(x_1, \ldots, x_N)$ are too complex to be solved algebraically and $f$ cannot be expressed as an algebraic function of $x_k$.

## 4 Examples

The examples we propose are simple in order to explain the resolution method without too long calculations. The Taylor series of $f$ can be calculated in any point, but they will be given only around the optimum points. Indeed, it is often of interest to evaluate the influence of the variables around the optima to get an idea of the accuracy required to reach them.

### 4.1 Example (1a)

Let us use the example of Fig. 1a to compare our method with the usual Lagrange multiplier method. We want to optimize $f(x, y) = x^2 + 2y^2$ with the constraint $(x, y) = x - y^2 = 1$ .

The Lagrangian function is $\mathcal{L}(x,y,\lambda) = x^2 + 2y^2 + \lambda(x - y^2 - 1)$. Its gradient is $\nabla\mathcal{L}$ is null when
$$\begin{cases} 2x + \lambda = 0 \\ 2y(2 - \lambda) = 0. \\ x - y^2 = 1 \end{cases}$$
This system is of dimension 3. There are three solutions
$$\begin{pmatrix} x = 1 \\ y = 0 \\ \lambda = -1/2 \end{pmatrix}, \begin{pmatrix} x = -1 \\ y = i\sqrt{2} \\ \lambda = 2 \end{pmatrix}, \begin{pmatrix} x = -1 \\ y = -i\sqrt{2} \\ \lambda = 2 \end{pmatrix}.$$
Only the first one is real and appears in Fig. 1a.

With the method we propose, the Jacobian is $J_{(f,g)} = \begin{pmatrix} 2x & 4y \\ 1 & -2y \end{pmatrix}$. The equation $Det\, J_{(f,g)} = 0$ makes the optimization problem of dimension 2 a solving problem of dimension 2:
$$\begin{cases} -4xy - 4y = 0 \\ x - y^2 = 1 \end{cases}.$$
The three solutions are $A = \begin{pmatrix} 1 \\ 0 \end{pmatrix}$, $B = \begin{pmatrix} -1 \\ i\sqrt{2} \end{pmatrix}$, $C = \begin{pmatrix} -1 \\ -i\sqrt{2} \end{pmatrix}$.

Even if it is here trivial, let us illustrate the role of the constraint matrices to determine the nature of the optimum in (1,0). Here, the Jacobian of the constraint is $J_{(g)} = [1, -2y]$, so its square submatrices are $S_1(x,y) = [-2y]$ and $S_2(x,y) = [1]$. The infinitesimal entities $dx$ and $dy$ are thus linked by $dy = \frac{-1}{-2y}dx$ and equivalently $dx = -(-2y)dy$.

The equation $dg = dx + 2ydy$ gives the boundaries of $\mathcal{C}_{(g)}$ the constraint domain curve. The $x$-boundary points are obtained by $dx = 0$ and $dy \neq 0$, which imposes that $y = 0$; the point is thus (1,0). This point is also the point of optimality of $f$. The $y$-boundary points are obtained by $dy = 0$ and $dx \neq 0$, which is impossible, and means that the curve is not bounded along the $y$-axis. This result is obvious is one considers the parabolic shape formed by the points of $g(x,y) = x - y^2 = 1$.

For the $x$-direction, $df = 2xdx + 4ydy = 2xdx + 2dx = 2(x+1)dx$. Thus, $f_x' = \frac{df}{dx} = 2(x+1)$, and $f_x'' = 2$. There is no other higher order terms $\frac{d^{(n)}f}{dx^n}$. Consequently, at the real optimal point (1,0), $f$ is exactly
$$\tilde{f}(x) = f(1,0) + f_x'(1,0).(x-1) + \frac{1}{2}f_x''(1,0).(x-1)^2 = 1 + 4(x-1) + (x-1)^2 = x^2 + 2x - 2$$
which is exactly the equation that could be obtained by substituting the expression of $y$ extracted from $g(x,y)$ inside the expression of $f(x,y)$.

For the $y$-direction, $df = 2xdx + 4ydy = 4xydy + 4ydy$. Thus, $f_y' = \frac{df}{dy} = 4y(x+1)$. With the partial derivative $df_y' = 4ydx + 4(x+1)dy = (8y^2 + 4x + 4)dy$. Thus, $f_y'' = 8y^2 + 4x + 4$. We conclude that at the optimal point (1,0), $f$ can be approximated by
$$\tilde{f}(y) \approx f(1,0) + f_y'(1,0).(y-0) + \frac{1}{2}f_y''(1,0).(y-0)^2 = 1 + 4y^2$$

The same result could be obtained by writing the Taylor series of $f(y) = 2y^2 + (1+y^2)^2$ up to the second order. The main advantage of the method we propose is that it works even when it is not possible to substitute the variables to write $f(x,y)$ as a function of $x$ uniquely, or as a function of $y$ uniquely, which can occur if the constraint equation contains complex intricate terms, for example $xy^5$ terms.

## 4.2 Example 1b

In the example of Fig. 1b we look for the optima of $f(x,y) = x^2 + 2y^2 - 2xy^2$ with the constraint $g(x,y) = x + y^2 = 1$

The Jacobian of the problem is $J_{(f,g)} = \begin{bmatrix} 2x - 2y^2 & 4y - 4xy \\ 1 & 2y \end{bmatrix}$. With the equation $Det\, J_{(f,g)} = 0$, the optimum problem becomes the solving problem

$$\begin{cases} 4xy - 4y^3 - 4y + 4xy = 0 \\ x + y^2 = 1 \end{cases}$$

The three solutions are $A = \begin{pmatrix} 1 \\ 0 \end{pmatrix}, B = \begin{pmatrix} 2/3 \\ 1/\sqrt{3} \end{pmatrix}, C = \begin{pmatrix} 2/3 \\ -1/\sqrt{3} \end{pmatrix}$.

To determine the nature of the optima, we use the Jacobian of the constraint $J_{(g)} = [1, 2y]$. Its submatrices are $S_1(x,y) = [2y]$ and $S_2(x,y) = [1]$. The infinitesimal enties $dx$ and $dy$ are thus linked by $dy = \frac{-1}{2y} dx$ and equivalently $dx = (-2y)dy$. Of course this relation could have been found directly by derivation of $g(x,y)$. The boundaries of the constraint domain $\mathcal{C}_{(g)}$ are the same as in the previous example and reduce to a unique point $A = \begin{pmatrix} 1 \\ 0 \end{pmatrix}$.

For the $x$-direction, $df = (2x - 2y^2)dx + (4y - 4xy)dy = (2x - 2y^2)dx + (2x - 2)dx$. Thus, $f'_x = \frac{df}{dx} = 4x - 2y^2 - 2$, and $df'_x = 4dx - 4ydy = 6dx$, i.e. $f''_x = 6$. There is no other higher terms. Consequently, $f$ is exactly

$$\tilde{f}(x) = f(x^P, y^P) + f'_x(x^P, y^P) \cdot (x - x^P) + \frac{1}{2} f''_x(x^P, y^P) \cdot (x - x^P)^2$$

For the $y$-direction, $df = (2x - 2y^2)dx + (4y - 4xy)dy = 2y(2y^2 - 2x)dy + (4y - 4xy)dy$. Thus, $f'_y = \frac{df}{dy} = 4y^3 - 8xy + 4y$. We continue to use the partial derivative

$df'_y = -8ydx + 4(3y^2 - 8x + 4)dy = (28y^2 - 8x + 4)dy$. Thus, $f''_y = 28y^2 - 8x + 4$. The polynomial approximation of $f$ around the three optimum points $P = (x^P, y^P)$ is

$$\tilde{f}(y) = f(x^P, y^P) + f'_y(x^P, y^P) \cdot (y - y^P) + \frac{1}{2} f''_y(x^P, y^P) \cdot (y - y^P)^2$$

For $A = \begin{pmatrix} 1 \\ 0 \end{pmatrix}$

$$\tilde{f}(x) = 1 + 2(x - 1) + 3(x - 1)^2 = 3x^2 - 4x + 2$$
$$\tilde{f}(y) = 1 - 2y^2$$

For $B = \begin{pmatrix} 2/3 \\ 1/\sqrt{3} \end{pmatrix}$

$$\tilde{f}(x) = \frac{2}{3} + 3\left(x - \frac{2}{3}\right)^2 = 3x^2 - 4x + 2$$

$$\tilde{f}(y) = \frac{2}{3} + 4\left(y - \frac{1}{\sqrt{3}}\right)^2 = 2 - \frac{8y}{\sqrt{3}} + 4y^2$$

For $C = \begin{pmatrix} 2/3 \\ -1/\sqrt{3} \end{pmatrix}$

$$\tilde{f}(x) = \frac{2}{3} + 3\left(x - \frac{2}{3}\right)^2 = 3x^2 - 4x + 2$$

$$\tilde{f}(y) = \frac{2}{3} + 4\left(y + \frac{1}{\sqrt{3}}\right)^2 = 2 + \frac{8y}{\sqrt{3}} + 4y^2$$

The parabolic surfaces of $\tilde{f}(x)$ and $\tilde{f}(y)$ are represented in Fig. 3 for the optimum points A and B. Those of C are deduced by symmetry $y \leftrightarrow -y$.

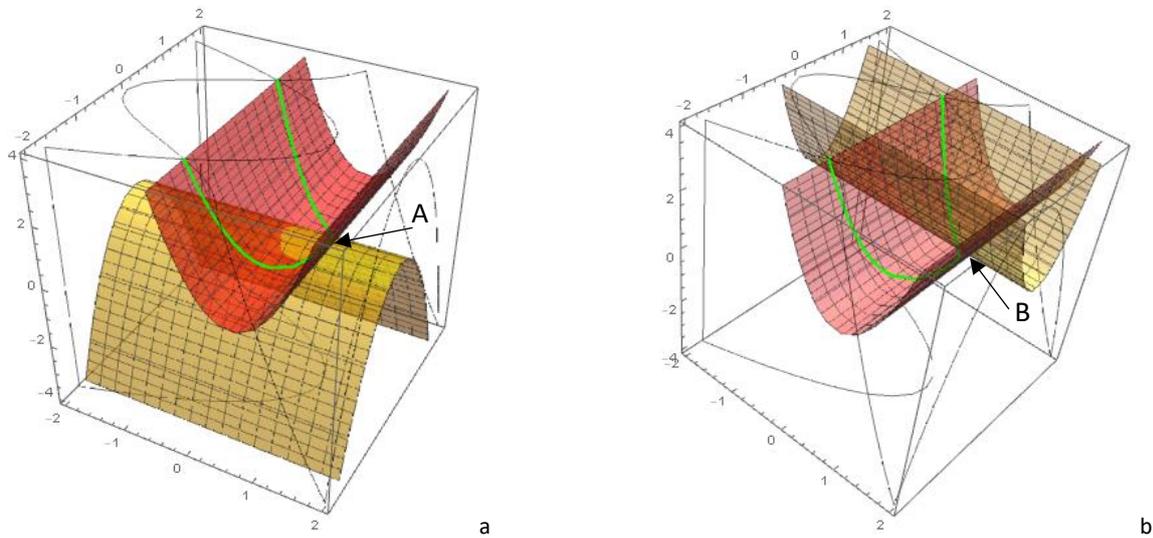

Fig. 3. Polynomial Taylor approximation functions $\tilde{f}(x)$ and $\tilde{f}(y)$ in red and yellow, respectively, of the function $f$ around its g-constraint optimum points, (a) A and (b) B.

### 4.3 Example 2

Let us propose now an example in 3D. We look for the optima of $f(x, y, z) = x^2 - 2y + z^3$ with the two constraints $g_1(x, y, z) = x^2 + y + z = 1$ and $g_2(x, y, z) = y - z^2 = -1$.

The contour surfaces of $f(x, y, z) = c$ are illustrated for different values of $c$ in Fig. 4. The contour surface of $g_1(x, y, z) = 1$ and that of $g_2(x, y, z) = -1$, with their intersection are shown in Fig. 5.

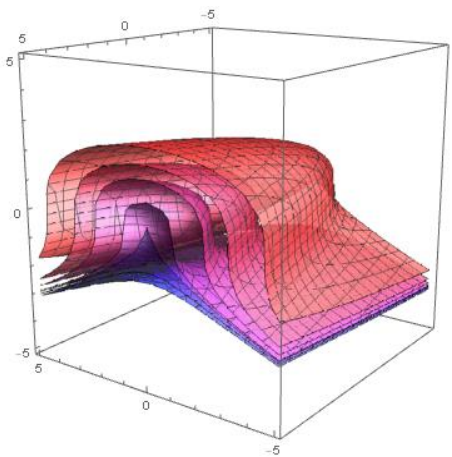 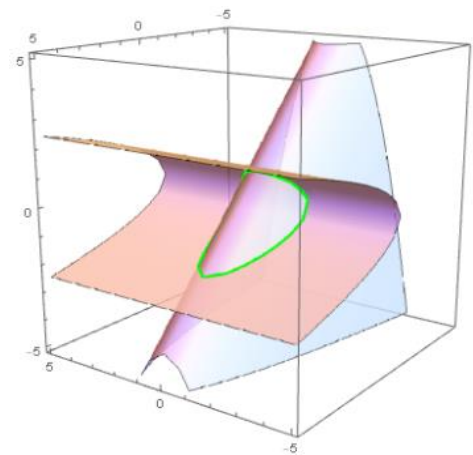

Fig. 4. Contour surfaces of $f(x, y, z) = x^2 - 2y + z^3 = c$ for $c \in \{-14, -12, -10, -9, -5, 0, 10\}$ from blue to red.

Fig. 5. Coutour surfaces of $g_1(x, y, z) = x^2 + y + z = 1$ and $g_2(x, y, z) = y - z^2 = -1$, with the intersection curve enlightened in green.

The Jacobian of the problem is $J_{(f,g_1,g_2)} = \begin{bmatrix} 2x & -2 & 3z^2 \\ 2x & 1 & 1 \\ 0 & 1 & -2z \end{bmatrix}$. The equation $Det J_{(f,g_1,g_2)} = 0$ makes the optimization problem of dimension 3 a solving problem of dimension 3:

$$\begin{cases} 6x(z-2)z - 2x = 0 \\ x^2 + y + z - 1 = 0 \\ y - z^2 + 1 = 0 \end{cases}$$

In general, the solution is found by numerical methods, but here the problem can be solved algebraically (we used Mathematica). There are four real solutions and two irrational ones. The four real optima $P_i$ are

$A = (0,3,-2)$, $B = (0,0,1)$, $C = \left(-\sqrt{2\sqrt{3}-\frac{4}{3}}, \frac{4}{3}(1-\sqrt{3}), 1-\frac{2}{\sqrt{3}}\right)$, $D = \left(\sqrt{2\sqrt{3}-\frac{4}{3}}, \frac{4}{3}(1-\sqrt{3}), 1-\frac{2}{\sqrt{3}}\right)$

The optimum points $P_i$ are positioned at the intersection of the three surfaces $f(x,y,z) = f(P_i)$ with the two constraints $g_1(x,y,z) = g_1(P_i) = 1$ and $g_2(x,y,z) = g_2(P_i) = -1$, as shown in Fig. 6.

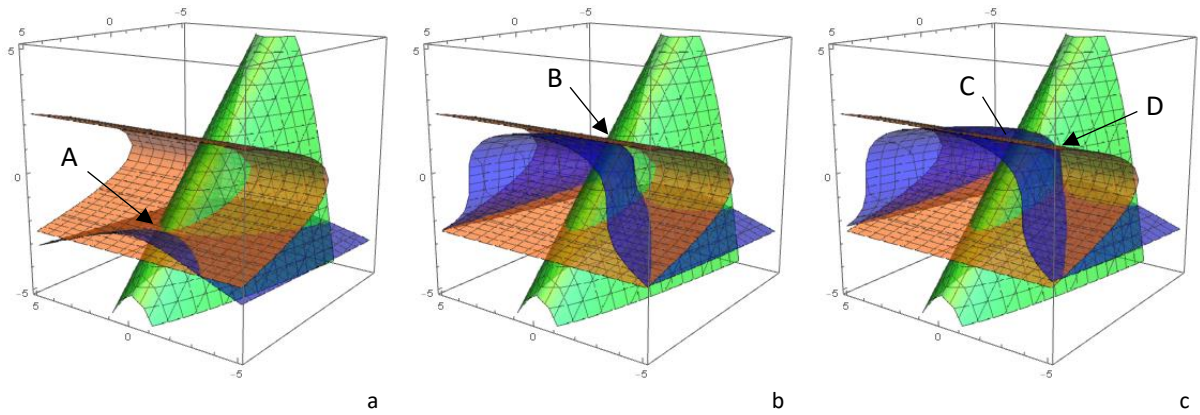

Fig. 6. Contour surfaces the two constraints $g_1(x,y,z) = g_1(P_i) = 1$ (in green) and $g_2(x,y,z) = g_2(P_i) = -1$ (in orange), with $f(x,y,z) = f(P_i)$ (in blue) at the optimum points $P_i =$ (a) A, (b) B, and (c) C and D. The values of $f$ at these points are (a) -14, (b) 1, (c) $1 + \frac{16}{3\sqrt{3}}$.

The linear relations between the infinitesimals $dx, dy, dz$ are given by the submatrices of $J_{(g_1,g_2)}$

$$J_{(g_1,g_2)} = \begin{bmatrix} 2x & 1 & 1 \\ 0 & 1 & -2z \end{bmatrix}$$

which are $S_x = \begin{bmatrix} 1 & 1 \\ 1 & -2z \end{bmatrix}$, $S_y = \begin{bmatrix} 2x & 1 \\ 0 & -2z \end{bmatrix}$, $S_z = \begin{bmatrix} 2x & 1 \\ 0 & 1 \end{bmatrix}$.

The boundary points of the constraint domain $C_{(g_1,g_2)}$ along $x, y$ and $z$ are determined by solving the initial constraints conditions of $g_1$ and $g_2$ with the additional condition that $Det\, S_x = 0$, $Det\, S_y = 0$, or $Det\, S_z = 0$, respectively, i.e.

| Along x: | Along y: | Along z: |
|---|---|---|
| $\begin{cases} -1 - 2z = 0 \\ x^2 + y + z - 1 = 0 \\ y - z^2 + 1 = 0 \end{cases}$ | $\begin{cases} -4xz = 0 \\ x^2 + y + z - 1 = 0 \\ y - z^2 + 1 = 0 \end{cases}$ | $\begin{cases} 2x = 0 \\ x^2 + y + z - 1 = 0 \\ y - z^2 + 1 = 0 \end{cases}$ |
| $\rightarrow \left\{\left(\frac{3}{2}, -\frac{3}{4}, -\frac{1}{2}\right), \left(-\frac{3}{2}, -\frac{3}{4}, -\frac{1}{2}\right)\right\}$ | $\rightarrow \{(0,3,-2), (0,0,1), (-\sqrt{2},-1,0), (\sqrt{2},-1,0)\}$ | $\rightarrow \{(0,3,-2), (0,0,1)\}$ |

As two points in the solution list are repeated twice, the constraint domain $\mathcal{C}_{(g)}$ is bounded by the 6 points represented in Fig. 7.

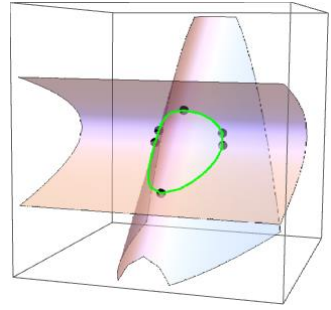

Fig. 7. The six points that bound the contraints domain.

The constraint matrices also allow us to create the linear links between the infinitesimals:

$$\begin{pmatrix}dy\\dz\end{pmatrix} = -S_x^{-1}\begin{pmatrix}2x\\0\end{pmatrix}dx \Rightarrow \begin{cases}dy = -\frac{4xz}{1+2z}dx\\dz = -\frac{2x}{1+2z}dx\end{cases} \qquad (\mathcal{R}_x)$$

$$\begin{pmatrix}dx\\dz\end{pmatrix} = -S_y^{-1}\begin{pmatrix}1\\1\end{pmatrix}dy \Rightarrow \begin{cases}dx = -\left(\frac{1}{2x}+\frac{1}{4xz}\right)dy\\dz = \frac{1}{2z}dy\end{cases} \qquad (\mathcal{R}_y)$$

$$\begin{pmatrix}dx\\dy\end{pmatrix} = -S_z^{-1}\begin{pmatrix}1\\-2z\end{pmatrix}dx \Rightarrow \begin{cases}dx = -\left(\frac{1}{2x}+\frac{z}{x}\right)dz\\dy = 2zdz\end{cases} \qquad (\mathcal{R}_z)$$

These relationships can now be used to determine the Taylor series of $f$ around its optimum points A, B, C, D. The derivative of $f$ is determined from the partial derivatives

$$df(x,y,z,dx,dy,dz) = 2xdx - 2dy + 3z^2 dz$$

Along the $x$-axis, thanks to the relation $(\mathcal{R}_x)$, we get

$df(x,y,z,dx,dy,dz) = 2xdx + 2\frac{4xz}{1+2z}dx - 3z^2\frac{2x}{1+2z}dx$. Thus, $f_x' = \frac{df}{dx}(x,y,z) = \frac{2x(1-3(-2+z)z)}{1+2z}$

Repeating the derivative with the partial derivative formula, we get

$$d(f_x') = \frac{2-6(-2+z)z}{1+2z}dx + x(-3+\frac{11}{(1+2z)^2})dz$$

Using again the relation $(\mathcal{R}_x)$, we get

$$d(f_x') = \frac{2-6(-2+z)z}{1+2z}dx + \left(-3+\frac{11}{(1+2z)^2}\right)\frac{-2x}{1+2z}dx$$

$$\Rightarrow f_x''(x,y,z) = \frac{-2(1+2z)^2(-1+3(-2+z)z)+8x^2(-2+3z(1+z))}{(1+2z)^3}$$

Along the $y$-axis, thanks to the relation $(\mathcal{R}_y)$, with the same method, we get

$$f_y' = \frac{df}{dy}(x,y,z) = \frac{3z}{2} - \frac{1}{2z} - 3$$

$$\Rightarrow f_y''(x,y,z) = \frac{3z^2+1}{4z^3}$$

Along the $z$-axis, thanks to the relation $(\mathcal{R}_z)$, with the same method, we get

$$f_z' = \frac{df}{dz}(x,y,z) = 3(z-2)z - 1$$

$$\Rightarrow f_z''(x,y,z) = 6(z-1)$$

The Taylor series along the $x$-axis is (the details are skipped):

For $A = (0, 3, -2)$, $\tilde{f}(x) = -14 + \frac{23x^2}{3}$

For $B = (0, 0, 1)$, $\tilde{f}(x) = 1 + \frac{4x^2}{3}$

For $C = \left(-\sqrt{2\sqrt{3} - \frac{4}{3}}, \frac{4}{3}(1 - \sqrt{3}), 1 - \frac{2}{\sqrt{3}}\right)$, $\tilde{f}(x) = \frac{-18089 + 11100\sqrt{3} + 48x(2(-54 + 35\sqrt{3})\sqrt{-4 + 6\sqrt{3}} + (105 - 54\sqrt{3})x)}{(-9 + 4\sqrt{3})^3}$

For $D = \left(\sqrt{2\sqrt{3} - \frac{4}{3}}, \frac{4}{3}(1 - \sqrt{3}), 1 - \frac{2}{\sqrt{3}}\right)$, $\tilde{f}(x) = \frac{-18089 + 11100\sqrt{3} + 48x(2(54 - 35\sqrt{3})\sqrt{-4 + 6\sqrt{3}} + (105 - 54\sqrt{3})x)}{(-9 + 4\sqrt{3})^3}$

The graph of $f(x)$ and its Taylor series around these points are shown in Fig. 8.

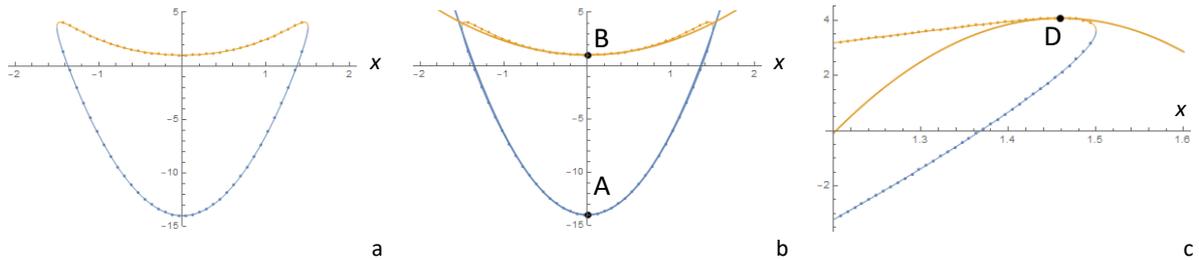

Fig. 8. Taylor series of $f$ along the x-axis. (a) Graph of $f(x)$. Here, the function $f(x, y, z)$ with the two contraints $g_1(x, y, z)$ and $g_2(x, y, z)$ can be transformed into two functions $f(x)$ by substitution, in blue and orange. (b) Parabolic Taylor series of $f(x)$ around A and B, and (c) around D (C is the same by symmetry $x \leftrightarrow -x$).

The Taylor series along the $y$-axis is:

For $A$, $\tilde{f}(y) = \frac{91}{64} - \frac{145y}{32} - \frac{13y^2}{64}$

For $B$, $\tilde{f}(y) = 1 - 2y + \frac{y^2}{2}$

For $C$ and $D$, $\tilde{f}(y) = -\frac{5}{9}(27 + 16\sqrt{3}) - \frac{1}{2}y(72 + 40\sqrt{3} + 36y + 21\sqrt{3}y)$

The graph of $f(y)$ and its Taylor series around these points are shown in Fig. 9.

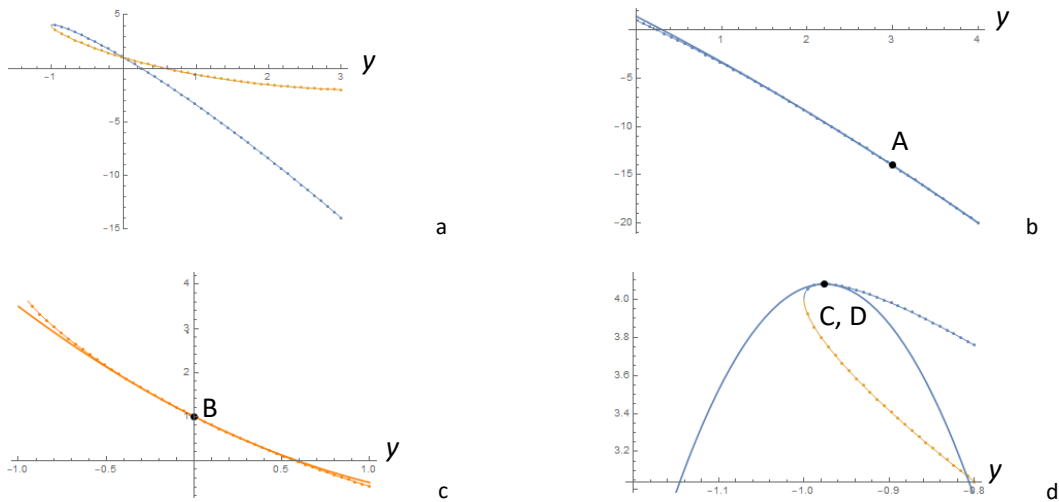

Fig. 9. Taylor series of $f$ along the $y$-axis. (a) Graph of $f(y)$. Here, the function $f(x, y, z)$ with the two contraints $g_1(x, y, z)$ and $g_2(x, y, z)$ can be transformed into two functions $f(y)$ by substitution, in blue and orange. (b) Parabolic Taylor series of $f(y)$ around A, (c) around B, and (d) around C and D.

The Taylor series along the $z$-axis is:

For $A = (0, 3, -2)$, $\tilde{f}(z) = \frac{499}{16} + \frac{355z}{16} - \frac{13z^2}{64}$

For $B = (0, 0, 1)$, $\tilde{f}(z) = \frac{11}{2} - 5z + \frac{z^2}{2}$

For $C = \left(-\sqrt{2\sqrt{3} - \frac{4}{3}}, \frac{4}{3}(1 - \sqrt{3}), 1 - \frac{2}{\sqrt{3}}\right)$, $\tilde{f}(z) = 1 + \frac{23}{6\sqrt{3}} - \frac{3}{2}z(4 + 2\sqrt{3} + 12z + 7\sqrt{3}z)$

For $D = \left(\sqrt{2\sqrt{3} - \frac{4}{3}}, \frac{4}{3}(1 - \sqrt{3}), 1 - \frac{2}{\sqrt{3}}\right)$, $\tilde{f}(x) = \frac{-18089 + 11100\sqrt{3} + 48x(2(54 - 35\sqrt{3})\sqrt{-4 + 6\sqrt{3}} + (105 - 54\sqrt{3})x)}{(-9 + 4\sqrt{3})^3}$

The graph of $f(z)$ and its Taylor series around these points are shown in Fig. 10.

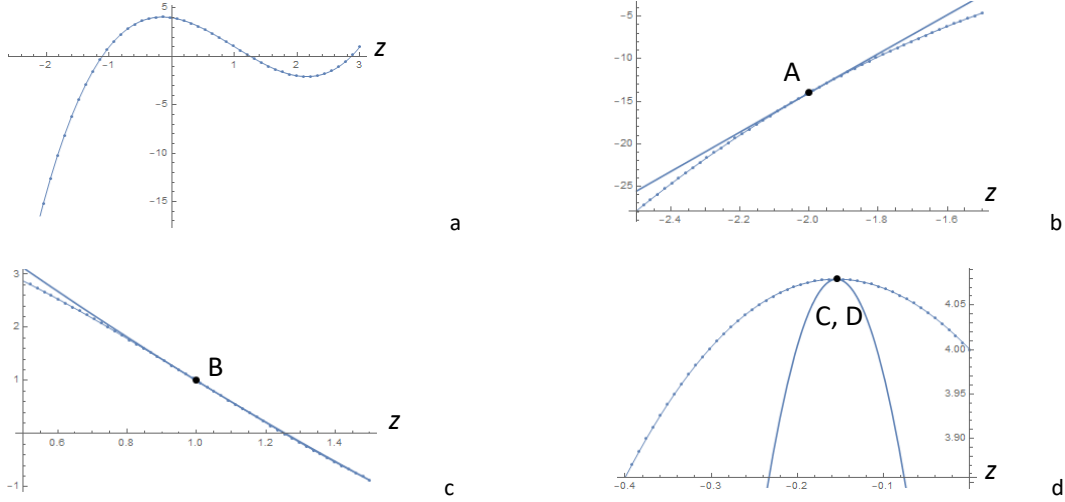

Fig. 10. *Taylor series of $f$ along the $z$-axis. (a) Graph of $f(y)$. Here, the function $f(x, y, z)$ with the two contraints $g_1(x, y, z)$ and $g_2(x, y, z)$ can be transformed into one function $f(z)$ by substitution. (b) Parabolic Taylor series of $f(z)$ around A, (c) around B, and (d) around C and D.*

## 5   Conclusion

The problem that consists in optimizing a $N$-dimension function with $M$ constraint equations is classically treated with the Lagrange multipliers and the associated Lagrange function. Generally, the problem is solvable if $M = N - 1$, and finding the optima is equivalent to solve a system of $2N - 1$ equations.

Here, we show that the problem that consists in optimizing a function $f$ of $N$ variables $x_1, \ldots, x_N$ with $N$-1 constraints of type $g_i(x_1, \ldots, x_N) = C_i$, $i \in [1, N - 1]$ can be solved without the Lagrange multipliers and without increasing the dimensionality of $f$. The optimality is obtained when $Det\, J_{(f,g)} = 0$, where $J_{(f,g)}$ is the Jacobian of the vectorial function $(f, g_1, \ldots, g_{N-1})$. The optimization problem is thus directly transformed into a solving problem of $N$ equations:

$$\begin{cases} \text{Optimize } f(x_1, \ldots, x_N) \\ g_1(x_1, \ldots, x_N) = C_1 \\ \vdots \\ g_{N-1}(x_1, \ldots, x_N) = C_{N-1} \end{cases} \Leftrightarrow \begin{cases} Det\, J_{(f,g)}(x_1, \ldots, x_N) = 0 \\ g_1(x_1, \ldots, x_N) = C_1 \\ \vdots \\ g_{N-1}(x_1, \ldots, x_N) = C_{N-1} \end{cases}$$

In addition, for any $k \in [1, N]$ we introduced $S_k$ the submatrices of the matrix $J_{(g)}$, the Jacobian of the constraint functions. They are obtained by removing the *k*-line and the *k*-column of $J_{(g)}$. The boundaries of the constraint curve $C_{(g)}$ along the $x_k$-axis are the points $(x_1, \ldots, x_N)$ that are the solutions of the system of $N-1$ equations:

$$\begin{cases} g_1(x_1, \ldots, x_N) = C_1 \\ \quad\vdots \\ g_{k-1}(x_1, \ldots, x_N) = C_{k-1} \\ Det\ S_k\ (x_1, \ldots, x_N) = 0 \\ g_{k+1}(x_1, \ldots, x_N) = C_{k+1} \\ \quad\vdots \\ g_{N-1}(x_1, \ldots, x_N) = C_{N-1} \end{cases}$$

Interestingly, the constraint matrices $S_k$ permit to establish linear relations between the infinitesimals $dx_i$. For any fixed $k \in [1, N]$, the infinitesimal $dx_i$ are proportional to $dx_k$ by

$$dx_i = s_{i,k}(x_1, \ldots, x_N)\, dx_k$$

with $s_{i,k}(x_1, \ldots, x_N) = -\left(S_k^{-1}\right)_{i,l} \frac{\partial g_l(x_1, \ldots, x_N)}{\partial x_k}$. Thanks to these relations and the partial derivative formula, one can determine by an iterative process for any $x_k$-axis the algebraic expression of $\frac{d^{(n)} f}{dx_k^n}(x_k)$, the *n*-derivative of $f$ along $x_k$. Consequently, for any point $(x_1, \ldots, x_N)_P$, $f$ can be developed in a Taylor series of $x_k$ around $(x_1, \ldots, x_N)_P$. The method works even if the algebraic forms of $f(x_1, \ldots, x_N)$ and $g_i(x_1, \ldots, x_N) = C_i$ are too complex to get the algebraic expression of $f(x_k)$ by substitution. Some simple 2D and 3D constrained optimization problems are proposed to show practically how to determine the optima of $f$, the boundaries of the constraint curve $C_{(g)}$, and the parabolic Taylor series of $f$ around the optimum points along the main $x_k$-axes.

We think that the method will simplify the resolution of some equality-constrained optimization problems. We would like to investigate its possible implications in physics, notably to clarify the connection between the second principle of thermodynamics and the least action principle of mechanics.